\documentclass[10pt,a4paper]{article}
\usepackage{amsmath,amsthm,amssymb}
\usepackage{color}
\usepackage[english]{babel}
\theoremstyle{plain}%
\newtheorem{theorem}{Theorem}[]%
\newtheorem{corollary}[theorem]{Corollary}%
\newtheorem{lemma}[theorem]{Lemma}%
\newcommand{\opH}{{H}}
\newcommand{\oph}{{h}}

\newcommand{\bigpar}{\par\quad\newline\noindent}

\newcommand{\norm}[1]{{\|{#1}\|}}
\newcommand{\abs}[1]{\left|{#1}\right|}

\newcommand{\ol}[1]{{\overline{#1}}}

\newcommand{\Rset}{{\mathbb{R}}}

\newcommand{\Nset}{{\mathbb{N}}}

\newcommand{\cointerval}[2]{[#1,\,#2)}%
\newcommand{\ccinterval}[2]{[#1,\,#2]}%

\newcommand{\at}[1]{{\left({#1}\right)}}

\newcommand{\bat}[1]{{\big(#1\big)}}

\newcommand{\ato}[1]{{\left[{#1}\right]}}

\newcommand{\pairo}[2]{{\left[{#1},\,{#2}\right]}}

\newcommand{\si}{{\sigma}}

\newcommand{\calM}{\mathcal{M}}

\newcommand{\calS}{\mathcal{S}}

\begin{document}
%

\title{Self-similar solutions with fat tails for a \\
coagulation equation with nonlocal drift}%
\date{\today}%
\author{%
Michael Herrmann\thanks{ %
Oxford Centre for Nonlinear PDE (OxPDE),
michael.herrmann@maths.ox.ac.uk}
\and
Philippe Lauren\c{c}ot\thanks{ %
Institut de Math\'ematiques de Toulouse,
laurenco@math.univ-toulouse.fr}
\and
Barbara Niethammer\thanks{ %
Oxford Centre for Nonlinear PDE (OxPDE),
niethammer@maths.ox.ac.uk}
}%
%
%
%
%
%
%
\maketitle
%
%
\begin{abstract}
We investigate the existence of self-similar solutions for a
coagulation equation with nonlocal drift. In addition to explicitly given
exponentially decaying solutions we establish the existence of
self-similar profiles with algebraic decay.
\end{abstract}
%
%
%
%
%
%
\section{Introduction}
%
%
The classical mean-field theory by Lifshitz \& Slyozov \cite{LS1} and
Wagner \cite{W1} describes the coarsening of droplets in dilute binary
mixtures and is based on the assumption that droplets interact only
via a common mean-field $u=u(t)$. It results in a nonlocal transport
equation for the number density $f=f(t,x)\ge 0$ of droplets of volume
$x>0$ at time $t\ge 0$ and reads  
\begin{equation}
\label{model2}
\partial_t f (t,x)+ \partial_x \big( (a(x) u(t) - b(x)) f(t,x) \big)=0,
\quad \qquad \int_0^{\infty} x f(t,x)\,dx = m_1 = \mathrm{const},
\end{equation}
where the second equation describes the conservation of matter
(volume) and determines the mean-field $u$. The functions $a$ and $b$
are specified by the mechanism of transfer of matter between droplets
and a typical example is $a(x)=x^\alpha$ and $b(x)=x^\beta$ with
$0\le\beta<\alpha\le 1$ (the original choice in \cite{LS1}
corresponding to $(\alpha,\beta)=(1/3,0)$).  
\par%
The large time behaviour of solutions to \eqref{model2} was
conjectured to be given by self-similar solutions already in
\cite{LS1,W1} but it is now known that it actually depends sensitively
on the details of the initial data at the end of their support
\cite{Carr1,CP1,NP2} . As a regularisation it was suggested in
\cite{LS1} to add a coagulation term with additive kernel to the
evolution equation for $f$ which accounts for the occasional merging
of droplets, that is,  
\begin{eqnarray}
\partial_t f(t,x) & + & \partial_x \big(  (a(x) u(t) - b(x)) f(t,x) \big) \nonumber\\
& = & \frac{1}{2}\ \int_0^x w(x-y,y) f(t,x-y) f(t,y)\,dy \,-\, f(t,x)
\int_0^{\infty} w(x,y) f(t,y)\,dy, \label{model1g} 
\end{eqnarray}
the mean-field $u$ still being given by the conservation of
volume. Well-posedness for this  case is proven in \cite{Laur3} and
the existence of a fast decaying self-similar solution is established
in \cite{HNV1} (for $a(x)=x^{1/3}$, $b(x)=1$, and $w(x,y)=x+y$); a
full characterization of all self-similar solutions seems however
difficult, as well as the study of their stability. In fact, these
questions are still open for the coagulation equation ($a=b=0$ in
\eqref{model1g}), except for the so-called solvable kernels $w(x,y)=2$
and $w(x,y)=x+y$ \cite{MP06}. In particular, besides the existence of
exponentially decaying self-similar profiles, nothing is known in
general about self-similar solutions with algebraic decay (``fat
tails'').  
\par%
In order to develop methods to tackle these problems we consider here
the following simplified version of \eqref{model1g}  
\begin{equation}
\label{model1}%
\partial_t f(t,x) + \partial_x \big(  (xu(t)-1  ) f(t,x) \big) =
\int_0^x f(t,x-y) f(t,y)\,dy \,-\, 2 f(t,x) \int_0^{\infty} f(t,y)\,dy,
\end{equation}
corresponding to the choice $a(x) = x$, $b(x)=1$, and $w(x,y) = 2$ for
$(x,y)\in (0,\infty)^2$. The function $u$ is again specified by
$\int_0^{\infty} x f(t,x)\,dx=m_1$, which is equivalent to
$u(t)=m_0/m_1$ with $m_0:=\int_0^{\infty} f(t,x)\,dx$. Self-similar
solutions to \eqref{model1}  are given by $f(t,x)=t^{-2}\ F (x/t)$ and
$u(t)= v/t$, for some $v \in (0,\infty)$. Introducing $z=x/t$ we
obtain that $(F,v)$ solves 
\begin{equation}
\label{model5}%
-\big( z(1-v) + 1\big) F'(z)= \big( 2 -v - 2 m_0 \big) F(z) + \int_0^z
F(z-y)F(y)\,dy\,, \qquad z\in (0,\infty)\,, 
\end{equation}
where $m_0=\int_0^{\infty} F(z)\,dz$, and $v$ is such that for given
$m_1>0$ the solution $F$ satisfies 
\begin{equation}
\label{model6} \int_0^{\infty} z F(z)\,dz = m_1,
\end{equation}
so \eqref{model5} implies $v=m_0/m_1$.
For the following analysis it is convenient to use $m_0$ and $v$ as
parameters. It is then easily seen that for each $v \in (0,1)$ there
is an exponentially decaying solution  
\[ 
F_v(z)= m_0 v e^{-vz} \qquad \mbox{ with } \qquad m_0=1-v.
\]
Besides these self-similar profiles which decay exponentially fast at
infinity,  we establish the existence of self-similar solutions with
algebraic decay provided that the parameter $m_0$ is sufficiently
small. 
\begin{theorem}
\label{P.1} %
For each $v\in (0,1)$ there exists $\overline{m_0}\in (0,1-v)$ such
that \eqref{model5}--\eqref{model6} has a unique solution $F$ with
$\int_0^\infty F(z) dz=m_0$ and $m_1=m_0/v$ for all
$m_0\in(0,\overline{m_0}]$. This solution $F$ is nonnegative,
non-increasing, and satisfies $F(z) \sim c z^{-(2-v)/(1-v)}$ as $z \to
\infty$ for some $c>0$. 
\end{theorem} 
In fact, we  conjecture that for all $m_0 \in (0,1-v)$ there is a
unique solution to \eqref{model5}--\eqref{model6} as in
Theorem~\ref{P.1} and that there is no solution with $m_0>1-v$. We aim
to prove this conjecture in a future work by a continuation method
starting from the solutions provided by Theorem~\ref{P.1}. 
%
%
\section{Proof of the existence result}
%
%
Within this section we always suppose that $v\in (0,1)$ is
fixed. Furthermore we denote by $M_i(F):= \int_0^{\infty} y^i
F(y)\,dy$ with $i\in\Nset$ the $i$th moment of a nonnegative and
integrable function $F$. 

%
\bigpar{\bf{Basic properties of solutions}. }%
%
We start with deriving some elementary properties of solutions of 
\eqref{model5}-\eqref{model6}. Integrating \eqref{model5} we easily 
establish the following Lemma.
\begin{lemma}
\label{L.ss1} Let $F$ be a solution of  \eqref{model5} such that
$M_0(F)=m_0$. Then we have $F(0)=m_0 (1-m_0)$, and
$m_1=M_1(F)<\infty$ implies $v=m_0/m_1$.
\end{lemma}
This result in particular implies that once we have established the
existence of a self-similar solution, then uniqueness follows by
uniqueness of the corresponding initial-value problem. Moreover, we
also infer from Lemma~\ref{L.ss1} that positive solutions can only
exist for $m_0 \in (0,1)$, but for technical reasons, and since our
results below requires $m_0$ to be sufficiently small anyway, we
assume from now on that $m_0<v/2$.
\bigpar
Since our existence proof relies on a fixed point argument for $F$ we
need appropriate supersolutions for \eqref{model5}. To this end we
define $\alpha:=(2-v-2m_0)/(1-v)\in (2,\infty)$ and the function
$\overline{F}_{m_0}(z):=m_0\ (1+(1-v)z)^{-\alpha}$ for $z\ge 0$, which
is the solution to the ordinary differential equation 
\begin{equation}
\notag
- \bat{z\at{1-v} +1}\overline{F}^\prime_{m_0}(z) = \big(2-v-2m_0 \big)
\overline{F}_{m_0}(z),\qquad \overline{F}_{m_0}(0)=m_0.
\end{equation}
This function satisfies $\overline{F}_{m_0}(z) \sim cz^{-\alpha}$ as
$z \to \infty$, and thanks to $m_0<v/2$  we find
$M_1(\overline{F}_{m_0}) < \infty$ as well as
$\overline{F}^\prime_{m_0}\at{z}\leq0$ for all $z\geq{0}$. As a
consequence of the maximum principle for ordinary differential
equations we readily derive the following comparison result. 
\begin{lemma}
\label{L.ss2} Any nonnegative solution $F$ to \eqref{model5} with
$F\at{0}\leq{m_0}$ satisfies $F(z) \leq \overline{F}_{m_0}(z)$ for
all $z\geq0$.
\end{lemma}
>From now on we restrict our considerations to \emph{admissible}
functions $F\in\mathcal{A}$, where $\mathcal{A}$ is the set of all
nonnegative and continuous functions
$F:\cointerval{0}{\infty}\to\cointerval{0}{\infty}$ with finite
moments $M_0\at{F}<\infty$ and $M_1\at{F}<\infty$. 
%
\bigpar{\bf{An auxiliary problem}. }%
%
A key ingredient in our existence proof is to show the existence of
solutions to the following auxiliary problem: For a given
$G\in\mathcal{A}$ with $M_0\at{G}=m_0$ we seek $F\in\mathcal{A}$ such
that 
\begin{equation}
\label{ss3a}
- \big( (1-v)z + 1 \big) F'(z)= \big (2-v-2m_0\big) F(z) + 2
\int_0^{z/2} F(z-y) G(y)\,dy\,, \qquad M_0(F)=m_0\,.
\end{equation}
Notice that the convolution operator in \eqref{ss3a} is related to
the integration over the interval $\ccinterval{0}{z/2}$. For $F=G$,
however, the identity $2 \int_0^{z/2} F(z-y) F(y)\,dy=\int_0^{z}
F(z-y) F(y)\,dy$ implies the equivalence of \eqref{ss3a} and
\eqref{model5}.
\begin{lemma}
\label{Lem.PropsSol} Consider $G\in\mathcal{A}$ with $M_0(G)=m_0$ and
suppose that $F\in\mathcal{A}$ solves \eqref{ss3a}. Then $F$ is
unique, monotonically decreasing, and fulfils
$m_0-2 m_0^2\leq{F}\at{0}\leq{m_0}$. 
\end{lemma}
\begin{proof}
The uniqueness of solutions follows from the homogeneity of the
problem combined with a  Gronwall-like argument, and the
monotonicity is implied by $F\geq0$, $G\geq 0$, $v<1$, and $2-v-2m_0\geq0$.
To derive the bounds for the initial value we integrate over
$z\in(0,\infty)$ to obtain
\begin{align*}
F(0) + (1-v) m_0 = (2-v) m_0 - 2m_0^2 + M_0 ( F *G) \;\;\mbox{ with
}\;\; (F*G)(z):= 2 \int_0^{z/2} F(z-y) G(y)\,dy, 
\end{align*}
and the estimates
\begin{align*}
0\leq{}M_0( F*G) = 2 \int_0^{\infty} \int_0^{z/2} F(z-y)G(y)\,dy \,dz
=  2 \int_0^{\infty} \int_0^z G(y) F(z)\,dy\,dz \leq 2 m_0^2
\end{align*}
completes the proof.
\end{proof}
For the subsequent considerations it is convenient to introduce the
function $\tau:\cointerval{0}{\infty}\to\Rset$ defined by
\begin{align}
\label{Def.Tau}
\tau\at{z}:=-z\frac{F'(z)}{F(z)} \quad\text{ as long as }\quad F(z)>0, %
\end{align}
so that
\begin{align}
\label{Sol.Tau} F(z)=F(0) \exp\at{- \int_0^z \frac{\tau(s)}{s}\,ds}.
\end{align}
Notice in particular that $\tau\at{z}\to\tau_\infty\neq0$ as $z \to
\infty$ implies $F\at{z} \sim c z^{-\tau_\infty}$ as $z \to
\infty$. Rewriting \eqref{ss3a} in terms of $\tau$ yields 
\begin{equation}
\label{ss5}
\frac{(1-v)z+1}{z} \tau(z) = 2 -v - 2m_0 +\oph\pairo{G}{\tau}\at{z},
\qquad
\oph\pairo{G}{\tau}\at{z} := 2 \int_0^{z/2} G(y)
\exp\at{\int_{z-y}^z \frac{\tau(s)}{s}\,ds} \,dy,
\end{equation}
and this is equivalent to the fixed point equation
\begin{equation}
\label{ss8}
\tau={\opH}\pairo{G}{\tau},\qquad{\opH}\pairo{G}{\tau}(z) := \frac{z}{(1-v)z +1 }
\big( 2-v-2m_0 + \oph\pairo{G}{\tau}(z) \big).
\end{equation}
In the following Lemma we summarize some useful properties of the operator $\opH$.
\begin{lemma}
\label{Lemma.HProps} For each $G\in\mathcal{A}$ the operator
$\tau\mapsto\opH\pairo{G}{\tau}$ is well defined on the set
$C\at{\cointerval{0}{\infty}}$ and has the following properties:
\begin{enumerate}
\item
$0\leq\opH\pairo{G}{\tau}$ for all $\tau$,
\item
$\tau_1\leq\tau_2$ implies $\opH\pairo{G}{\tau_1}\leq{\opH}\pairo{G}{\tau_2}$,
\item
$\tau\at{z}\to\tau_\infty$ as $z\to\infty$ implies
${\opH}\pairo{G}{\tau}\at{z}\to (2-v)/(1-v)$ as  $z\to\infty$.
\end{enumerate}
\end{lemma}
\begin{proof}
The first two claims are direct consequences of $G \geq 0$ combined with
$v<1$ and $m_0\leq v/2<1/2$, and the definitions \eqref{ss5} and
\eqref{ss8}. Now suppose that 
$\tau\at{z}\to\tau_\infty$, and write
\begin{align*}
\oph\pairo{G}{\tau}\at{z} = 2 \int_0^{z/2} G(y) \exp \Big(
\int_{z-y}^z \frac{\tau(s)-\tau_{\infty}}{s}\,ds \Big) \Big(
\frac{z}{z-y} \Big)^{\tau_{\infty}}\,dy.
\end{align*}
Due to $y\leq{z}/2$ the convergence assumption implies
\begin{align*}
\abs{\exp \Big( \int_{z-y}^z \frac{\tau(s)-\tau_{\infty}}{s}\,ds
\Big)-1} \leq \exp \Big ( \ln 2
\sup\limits_{s\geq{z}/2}\abs{\tau\at{s} - \tau_{\infty}} \Big)-1
\quad\xrightarrow{\quad{z}\to \infty\quad}\quad0,
\end{align*}
and from $M_0\at{G}=m_0<\infty$, the inequality $z/(z-y)\le 2$ for
$y\in (0,z/2)$, and the Lebesgue Dominated Convergence Theorem we
infer that 
\begin{align*}
\lim_{z\to \infty} 2 \int_0^{z/2} G(y) \Big( \frac{z}{z-y} \Big
)^{\tau_{\infty}}\,dy = 2 \int_0^\infty G(y) \,dy  = 2m_0.
\end{align*}
Hence, we have $\oph\pairo{G}{\tau}\at{z}\to2m_0$ as $z\to\infty$ ,
and the proof is complete.
\end{proof}
%
%
%
\bigpar{\bf{Solutions to the auxiliary problem}. }%
%
In order to set up an iteration scheme for $\tau$ we prove that for
all sufficiently small $m_0$ there exists a supersolution to the fixed
point equation \eqref{ss8}. Notice that the constants in the next
Lemma depend on $v$.
\begin{lemma}
\label{Lemma.Upperbound}%
There exist two constants $\overline{m_0}\in(0,v/2)$ and $\tau_\star\in\Rset$
such that $\opH\pairo{G}{\tau_\star}\leq{\tau_\star}$ holds for all
$G\in\mathcal{A}$ with $m_0=M_0\at{G}\leq\overline{m_0}$.
\end{lemma}
\begin{proof}
The definition \eqref{ss5} gives
\begin{align}%
\notag
\oph\pairo{G}{\tau}\at{z} \leq 2 \int_0^{z/2} G(y) \Big(
\frac{z}{z-y} \Big)^{\|\tau\|_{\infty}}\,dy \leq 2 m_0
2^{\|\tau\|_{\infty}}
\end{align}
and hence $\opH\pairo{G}{\tau} \leq
a_{m_0}+b_{m_0}2^{\norm{\tau}_\infty}$ with $a_{m_0}:=(2 - v -
2m_0)/(1-v)$ and $b_{m_0}:=2m_0/(1-v)$. Since $b_{m_0}\to 0$ as
$m_0\to 0$ we choose $\ol{m_0}>0$ and $\si_\star>0$ such that 
$b_{\overline{m_0}}2^{a_0}\leq\si_\star2^{-\sigma_\star}$, and set
$\tau_\star=a_0+\sigma_\star$. Then, for all $m_0\in (0,\overline{m_0}]$
we find
\begin{math}
a_{m_0}+b_{m_0}2^{\tau_\star}
\leq
{a_0}+b_{\overline{m_0}}2^{a_0}2^{\si_\star}\leq{a_0}+{\sigma^\star}=\tau_\star,
\end{math} %
and this implies the claimed result.
\end{proof}
Now we are able to prove the existence of solutions to the auxiliary
problem.
\begin{corollary}
\label{P.2}
For each $G\in\mathcal{A}$ with $m_0=M_0(G)\in (0,\overline{m_0}]$ there
exists a unique solution $F\in\mathcal{A}$ to \eqref{ss3a} such
that the function $\tau$ given by \eqref{Def.Tau} satisfies
$0\leq\tau\at{z}\leq\tau_\star$ for all $z\geq0$ and
$\tau\at{z}\to(2-v)/(1-v)$ as $z\to\infty$. In addition, $F\leq\overline{F}_{m_0}$.
\end{corollary}
\begin{proof}
We define a sequence of continuous functions $\at{\tau_i}_i$ by
$\tau_0=\tau_\star$ and $\tau_{i+1}=\opH\pairo{G}{\tau_i}$, which
satisfies $0\leq\tau_{i+1}\leq{\tau}_i\leq{\tau_\star}$ thanks to
Lemma~\ref{Lemma.HProps} and Lemma~\ref{Lemma.Upperbound}. Consequently,
this sequence converges pointwisely to a function $\tau$, which
solves the fixed point equation \eqref{ss8}, and is hence
continuous. The function $F$ is then determined by \eqref{Sol.Tau} 
and the condition $M_0\at{F}=m_0$, the latter being meaningful as the
behaviour of $\tau$ as $z\to\infty$ guarantees the integrability of
$F$. Moreover, $F$ satisfies the first equation in \eqref{ss3a} by
construction and $F\leq\overline{F}_{m_0}$ due to
Lemma~\ref{Lem.PropsSol}, and this implies $M_1\at{F}<\infty$. 
\end{proof}
%
%
\bigpar{\bf{Fixed point argument for $F$}. }%
%
We finish the proof of Theorem~\ref{P.1} by applying Schauder's Fixed
Point Theorem. For $m_0\in (0,\overline{m_0}]$ let $\calM_{m_0}$
denote the set of all functions $G\in\mathcal{A}$ that satisfy
$M_0\at{G}=m_0$ and $G\leq\overline{F}_{m_0}$. In view of
Lemma~\ref{Lem.PropsSol} and Corollary~\ref{P.2}, we may define an
operator $\calS_{m_0}:\calM_{m_0}\to\calM_{m_0}$ as follows: For given
$G\in\calM_{m_0}$ the function $\calS_{m_0}\ato{G}$ is the solution to
the auxiliary problem \eqref{ss3a} with datum $G$. 
\begin{lemma}
For each $m_0\in (0,\overline{m_0}]$ the operator $\calS_{m_0}$ has a
unique fixed point in $\calM_{m_0}$, which satisfies both
\eqref{model5} and \eqref{model6} with $m_1=m_0/v$.
\end{lemma}
\begin{proof}
>From \eqref{ss3a} and Lemma~\ref{Lem.PropsSol} we infer that each
$F\in\calS\ato{\calM_{m_0}}$ is differentiable with 
\begin{align*}
\norm{F^\prime}_\infty\leq\at{2-v+2m_0}\norm{F}_\infty\leq\at{2-v+2m_0}{m_0}.
\end{align*}
Therefore, $\calS_{m_0}\ato{\calM_{m_0}}$ is precompact in
$C(\cointerval{0}{\infty})$, and due to the 
uniform supersolution $\overline{F}_{m_0}$ we readily verify that
the integral constraint $M_0\at{G}=m_0$ is preserved under strong
convergence in $\calM_{m_0}$. Schauder's theorem now implies the
existence of a fixed point $F=\calS_{m_0}\ato{F}\in{\calM}_{m_0}$,
which satisfies \eqref{model5} and $M_0(F)=m_0$ by construction. Moreover,
$F\leq\overline{F}_{m_0}$ implies $M_1\at{F}<\infty$, so
Lemma~\ref{L.ss1} guarantees \eqref{model6}. Finally, the fixed point
is unique as discussed in the remark to Lemma~\ref{L.ss1}. 
\end{proof}


\begin{thebibliography}{00}
%
\bibitem{Carr1}
J.~Carr.
\newblock Stability of self-similar solutions in a simplified {LSW} model.
\newblock {\em Phys. D}, 222(1-2):73--79, 2006.

\bibitem{CP1}
J.~Carr and O.~Penrose.
\newblock {Asymptotic behaviour in a simplified Lifshitz--Slyozov equation}.
\newblock {\em Phys. D}, 124:166--176, 1998.

\bibitem{HNV1}
M.~Herrmann, B.~Niethammer, and J.J.L. Vel\'azquez.
\newblock {Self-similar solutions for the LSW model with encounters}.
\newblock 2007.
\newblock preprint.

\bibitem{Laur3}
Ph.~Lauren{\c{c}}ot.
\newblock The {L}ifshitz-{S}lyozov equation with encounters.
\newblock {\em Math. Models Methods Appl. Sci.}, 11(4):731--748, 2001.

\bibitem{LS1}
I.~M. Lifshitz and V.~V. Slyozov.
\newblock The kinetics of precipitation from supersaturated solid solutions.
\newblock {\em J. Phys. Chem. Solids}, 19:35--50, 1961.

\bibitem{MP06}
G.~Menon and R.~L. Pego.
\newblock {Approach to self-similarity in Smoluchowski's coagulation
equations}.
\newblock {\em Comm. Pure Appl. Math.}, 57(9):1197--1232, 2004.

\bibitem{NP2}
B.~Niethammer and R.~L. Pego.
\newblock {Non--self--similar behavior in the LSW theory of Ostwald ripening}.
\newblock {\em J. Statist. Phys.}, 95(5/6):867--902, 1999.

\bibitem{W1}
C.~Wagner.
\newblock{Theorie der Alterung von Niederschl\"agen durch Uml\"osen
(Ostwald-Reifung)}. 
\newblock{\em Z. Elektrochem.}, 65:581--591, 1961.

\end{thebibliography}
\end{document}